\documentclass[10pt,a4paper]{article}
\usepackage[utf8x]{inputenc}
\usepackage[english]{babel}
\usepackage{lmodern}
\usepackage{amsmath}
\usepackage{amssymb}
\usepackage{amsthm}
\usepackage{enumerate}
\usepackage{subfigure}
\usepackage{bm}
\usepackage{anysize}
\usepackage{graphicx}
\usepackage{epsfig}
\usepackage{latexsym}
\usepackage{enumerate,cases}
\usepackage{multirow}
\usepackage{graphicx}
\usepackage{oxford3}
\newif\iffancyhd
\fancyhdtrue 

\newif\ifhdlowcase 

\iffancyhd

\def\i{{\rm i}}

\def\E{\mathrm E}
\def\P{\mathrm P}
\def\e{{\rm e}}
\def\d{{\rm {d}}}
\def\Z{\mathbb{Z}}

\def\R{\mathbb{R}}

\def\1{{\bf 1}}
\def\text{\mbox}

\def\eq2{
\stackrel{\small \rm mod \,2}{=}}
\def\n2{
\stackrel{\small \rm mod \,2}{\neq}}

\newcommand{\toD}{\stackrel{\rm  d}{\rightarrow}}
\newcommand{\toFDD}{\stackrel{\rm  fdd}{\rightarrow}}
\newcommand{\toP}{\stackrel{\rm P}{\rightarrow}}
\newcommand{\limfdd}{\stackrel{\rm  fdd}{\rightarrow}}

\newcommand{\eqfdd}{\stackrel{\rm  fdd}{=}}

\newcommand{\mbf}[1]{\mbox{\boldmath $#1$}}

\newcommand{\dd}{\,\mathrm{d}} 

\DeclareMathOperator{\Var}{Var}

\theoremstyle{plain}

\def\abstract#1{\vspace{2mm}{\raggedright\textbf{Abstract}. {#1}\vspace{1mm}}}

{\catcode `\@=11 \global\let\AddToReset=\@addtoreset}
\AddToReset{equation}{section}

\begin{document}
\title{Aggregation and long memory: recent developments\thanks{The first, third and fourth authors 
are supported by a grant (No.\ MIP-13079) from the Research Council of Lithuania.}
}
\author{R. Leipus$^{1}$, A. Philippe$^{2}$, D. Puplinskait\.e$^{1,2}$ and D. Surgailis$^{1}$}
\date{\today \\  \small
\vskip.2cm
$^1$Vilnius University and $^2$Universit\'{e}  de Nantes }
\maketitle

\begin{abstract}
It is well-known that the aggregated time series might have very different properties from those of the individual series, in particular,
long memory. At the present time, aggregation has become one of the main tools for modelling of long memory processes.
We review recent work on contemporaneous aggregation of random-coefficient AR(1) and related models,
with particular focus on various long memory properties of the aggregated process.  \\

\bigskip

\end{abstract}

\vskip1cm

\section{Introduction}
\label{sec:introduction}

In the seminal paper \citetext{bib:GR} observed that the covariance function of
AR(1) process with random and Beta distributed coefficient
may decay very slowly as in ARFIMA process. Indeed, let $X(t) = \sum_{j=0}^\infty a^j \zeta(t-j) $ be a stationary
solution of AR(1) equation
\begin{equation}\label{AR0}
X(t) = a X(t-1) + \zeta(t),\quad t\in \Z,
\end{equation}
where $\{ \zeta(t) \}  \sim WN(0,\sigma^2)$ is a white noise
and $a\in [0,1)$ is a r.v., independent of $\{\zeta(t) \}$, and having a density $\phi(x)$ regularly varying at the unit root $x=1$, viz.,
\begin{equation} \label{mix1}
\phi (x) \sim c_\phi (1-x)^\beta,  \quad x \nearrow 1,
\end{equation}
where $c_\phi >0 $ and $\beta >- 1$ some constants. Assume $0< \beta <1$ in the rest of this section. Then, as $t \to \infty$,
\begin{eqnarray}
\E X(0)X(t)&=&\sigma^2 \int_0^1 \frac{x^{t}}{1- x^2} \phi(x) {\d}x \
\sim \  c\,t^{-\beta}, \label{covX}
\end{eqnarray}
with $c= \sigma^2 (c_\phi/2) \Gamma(\beta)$ with $\Gamma(\beta) = \int_0^\infty {\e}^{-y} y^{\beta-1} {\d}y$. In the case
of Beta-type density  $\phi (x) = \frac{2}{B(p,q)}
x^{2p-1}(1-x^2)^{q-1}$ 
discussed in \citetext{bib:GR}, condition (\ref{mix1}) is satisfied with $\beta = q-1$ and the covariance
in (\ref{covX}) can be explicitly computed: $\E X(0)X(t)
= (\sigma^2 \Gamma(q-1)/B(p,q))\Gamma (p + \frac{t}{2})/\Gamma (p + \frac{t}{2} + q-1)$,
leading to the asymptotics in (\ref{covX}) since $\Gamma(t)/\Gamma(t+b) \sim t^{-b}$, $t \to \infty$.

Now suppose that one wants to study a huge and heterogeneous population of dynamic ``micro-agents" who all evolve independently of each other according to AR(1) process. Thus, the probability law of the  evolution of a given  ``micro-agent" is
completely determined by the value of the autoregressive parameter $a$. The heterogeneity of ``micro-agents" means
that $a$ has a probability density $\phi$ across the population.  The ``macroeconomic" variable $\mathfrak{X}_N(t)$  of interest is
obtained by averaging the evolutions of $N$ ``microeconomic" variables $X_i(t), i=1, \dots, N$,
\begin{equation}\label{XN}
\mathfrak{X}_N(t)\ = \ \frac{1}{A_N} \sum_{i=1}^N X_i(t), \quad  t =
0,\pm 1, \dots,
\end{equation}
which are randomly sampled from all ``micro-agents" population, where $A_N$ is a normalization and 
$\{ X_i(t)\}, i=1,2, \dots $ are i.i.d.\ copies of (\ref{AR0}).
By the independence of the summands
in (\ref{XN}) and the classical CLT, it immediately follows (with $A_N = N^{1/2}$) that
\begin{equation}\label{Xlim}
\mathfrak{X}_N(t) \ \toFDD \  \mathfrak{X}(t),
\end{equation}
where  $\{\mathfrak{X}(t)\}$ is a stationary Gaussian process with the same 2nd moment characteristics as the individual
``micro-agents", i.e.
\begin{equation}\label{covXG}
\E \mathfrak{X}(t)\ =\ 0, \qquad \E \mathfrak{X}(0) \mathfrak{X}(t)\ = \ \E X(0) X(t)
\ = \ \sigma^2 \int_0^1 \frac{x^{t}}{1- x^2} \phi(x) {\d}x
\end{equation}
and where $\toFDD$ denotes the weak convergence of finite-dimensional distributions. In particular, under
the assumption (\ref{mix1}), $0< \beta < 1$,
the covariance $\E \mathfrak{X}(0) \mathfrak{X}(t)$ decays as in (\ref{covX}), implying that
$\sum_{t \in \Z} |\E \mathfrak{X}(0) \mathfrak{X}(t)| = \infty,$ i.e., $\{\mathfrak{X}(t)\}$ has long memory.
Notice also that the Gaussian process
$\{\mathfrak{X}(t)\}$ in (\ref{Xlim}) is ergodic, in contrast to the non-ergodic random-coefficient individual  processes $\{X_i(t)\}$
in (\ref{XN}). Since (\ref{XN}) refers to summation of (independent) processes at instant time $t$,
the summation procedure (\ref{XN}) is called {\it contemporaneous} or
{\it cross-sectional }  aggregation, to be distinguished from {\it temporal aggregation},
the latter term usually referring to the summation on disjoint successive blocks, or taking partial sums, of a single time series.
We see that contemporaneous aggregation of simple dynamic models can result in a long memory process and hence may
provide an explanation of the long memory phenomenon observed in many econometric studies (\citetext{bib:baillie:1996},
\citetext{kirman-teyssiere}).
The last fact is very important since economic reasons of long memory remain largely unclear and some
authors try explain this phenomenon by ``spurious long memory'' (\citetext{Lobato:Savin-98}, \citetext{Mikosch-03}).

It is also instructive to look at the behavior (\ref{covX}) from the ``spectral perspective". The spectral density of
the random-coefficient AR(1) process is written as
$$
f(y)= \frac{\sigma^2}{2\pi}\int_0^1 \frac{1}{|1 - x\e^{\i y}|^2 } \phi(x) \d x
= \frac{\sigma^2}{2\pi}\int_0^1 \frac{1}{(1-x)^2 + 4x \sin^2 (y/2) } \phi(x) \d x, \quad y\in [-\pi,\pi].
$$
We see that under (\ref{mix1}), $0< \beta < 1$,
when $y\searrow  0$, this ``aggregated''  spectral density behaves like a power function:
\begin{eqnarray}
f(y)&\sim& \sigma^2\frac{c_\phi}{2\pi}\int_0^1 \frac{(1-x)^\beta \d x}{ (1- x)^2 + 4 x \sin^2 (y/2)}\
\sim\ \sigma^2 \frac{c_\phi}{2\pi}\int_0^1 \frac{x^\beta \, \d x}{x^2 + 4  \sin^2 (y/2)}\nonumber \\
&=&\sigma^2\frac{c_\phi} {2\pi (2 \sin (y/2))^{1-\beta}} \int_0^{1/2\sin(y/2)}   \frac{w^\beta \d w}{ w^2 + 1} \nonumber \\
&\sim& \frac{c_f}{y^{1-\beta}}, \label{specX}
\end{eqnarray}
where $c_f := \sigma^2 \frac{c_\phi}{2\pi} \int_0^{\infty}
\frac{w^\beta \d w}{ w^2 + 1}$. The singularity in the low frequency spectrum of
$\{X(t)\}$ and $\{\mathfrak{X}(t)\}$ is another indication of long memory of these processes and
is not surprising as the power behaviors (\ref{covX}) and
(\ref{specX}) are known to be roughly equivalent. Relation (\ref{specX}) seems to be well-known to physicists and
even experimentally observed, relating aggregation and long memory to the vast area in physics called
$1/f$ noise (see, e.g., the article by Ward and Greenwood
(\citeyear{ward-greenwood}) in Scholarpedia).
Another manifestation of long memory for the aggregated process in \eqref{specX} is the convergence of normalized partial
sums of $\{\mathfrak{X}(t)\}$ to a fractional Brownian motion (\citetext{bib:ZA}). See Sec.\ \ref{sec:aggr-ar1-proc} for the
related notion of distributional long memory.
\smallskip

Following \citetext{bib:GR}, various aspects of aggregation were discussed in the literature, see
\citetext{bib:GG}, Zaffaroni (\citeyear{bib:ZA}, \citeyear{bib:ZA07a}), \citetext{bib:OV}, \citetext{bib:LOPV06}, Celov {\em et al.} (\citeyear{bib:CE}, \citeyear{bib:CLP10}), \citetext{jirak:2011}. \citetext{bib:OV} considered aggregation of general
AR($2p)$ processes with $2p$ random coefficients, under the condition that the characteristic polynomial
has the form
 \begin{equation*}
 A(z)\ = \ (1-\alpha_{1} z)(1+\alpha_{2} z) \prod_{k=3}^{p+1} (1-\alpha_{k}{\rm e}^{{\rm i}\theta_k} z)
 (1-\alpha_{k}{\rm e}^{-{\rm i}\theta_k} z),
\end{equation*}
 where $\alpha_1,\dots,\alpha_{p+1}$ are i.i.d.\ random variables with probability densities as in (\ref{mix1}),
 with (possibly) different $\beta$'s, and $\theta_3,\dots,\theta_{p+1}$ are fixed. They described the singularities
of the spectral  density and the asymptotics of the covariance function of the corresponding aggregated Gaussian process, which
may contain a seasonal component. \citetext{bib:ZA} considered a general aggregation scheme of random-coeffcient AR(1)
with idiosyncratic and common components. Several papers (\citetext{bib:DG}, \citetext{bib:ZA07a},  \citetext{bib:ZA07b},
\citetext{bib:LV}, \citetext{bib:KLV04}, \citetext{bib:GLS08}) extended the aggregation scheme to ARCH-type  heteroskedastic processes
with random coeffcients and common innovations, with a particular emphasis on long memory behavior. The disaggregation problem
of reconstructing the mixing distribution from the spectral density or observed sample from the aggregated process
and its statistical aspects were studied \citetext{bib:DCF06}, \citetext{bib:CE}, \citetext{bib:LOPV06},  \citetext{bib:CLP10}, \citetext{jirak:2011},
\citetext{bib:ROB78}, \citetext{bib:BSG10} and elsewhere. Some aspects of aggregation and disaggregation are also discussed in
the recent monograph \citetext{bib:beran2013}.

\smallskip

The aim of the present paper is to review  recent developments on
aggregation and long memory. Sec.\ \ref{sec:aggr-ar1-proc} discusses the case of AR(1)
processes with infinite variance. Sec.\ \ref{sec:aggr-ar1-proc-2} extends the aggregation procedure
to the triangular array model. Sec.\ \ref{sec:joint} studies the joint temporal ant contemporaneous
aggregation of AR(1) processes. Sec.\ \ref{sec:fields} considers the aggregation procedure for
random fields. Sec.\ \ref{sec:disaggr} reviews some results related to statistical inference for the mixing density
$\phi$.

\section{Aggregation of  AR(1) processes with infinite variance}
\label{sec:aggr-ar1-proc}

The approach of \citetext{bib:GR} presented Sec.\ \ref{sec:introduction} can be extended to the case of infinite
variance. Let
\begin{equation}\label{AR1}
X(t) \ = \ \sum_{j=0}^\infty a^j \zeta(t-j) 
\end{equation}
be a stationary solution of the AR(1) equation (\ref{AR0})
with random coefficient $a \in [0,1)$ and
i.i.d.\ innovations  $\{ \zeta(t) \}$ with {\it infinite variance} ${\rm Var}(\zeta(t)) = \infty $, independent of r.v. $a$,
and belonging to the domain of attraction of $\alpha-$stable law, $0 < \alpha < 2$, in the sense that
\begin{equation}\label{DA}
\frac{1}{N^{1/\alpha}} \sum_{i=1}^N  \zeta(i)\ \toD\ Z,
\end{equation}
where $Z$ is an $\alpha-$stable r.v. Let $\{X_i(t)\}, i=1,2, \dots $ be independent copies of (\ref{AR1}), and define
the aggregated process as
\begin{equation}\label{XN1}
\mathfrak{X}_N(t)\ := \ \frac{1}{N^{1/\alpha}} \sum_{i=1}^N X_i(t), \quad  t\in \Z.
\end{equation}
The problem is to determine the limit process $\{\mathfrak{X}(t)\}$ of (\ref{XN1}), in the sense
of (\ref{Xlim}), and then we  describe its properties, in particular, long memory properties.

The above questions were studied in \citetext{bib:ps2010}.  Note that
the AR(1) series in (\ref{AR1}) converges conditionally a.s.\ for every $a \in [0,1)$, and also unconditionally a.s.\
if the distribution of $a$ satisfies the additional condition $\E (1-a)^{-1} < \infty. $ For regularly varying $\phi$ as in (\ref{mix1}),
the last condition is equivalent to $\beta > 0$. It is shown in Puplinskait\.e and Surgailis (\citeyear{bib:ps2010}, Thm.\ 2.1)
that in the case $\beta >0$ the limit of $\{\mathfrak{X}_N(t)\}$ exists and is written as a stochastic integral
\begin{equation}
\mathfrak{X}(t)\ :=\ \sum_{s \le t} \int_{[0,1)} x^{t-s} M_{\alpha,s}({\d}x),  \label{mixstable}
\end{equation}
where  $\{M_{\alpha,s}, s\in \Z\} $ are i.i.d. copies of an $\alpha-$stable random measure  on $[0,1)$ with distribution $Z$ in (\ref{DA})
and control measure
equal to the mixing distribution $\Phi $ (i.e., the distribution of the r.v.\ $a$). Recall that a family
$M =  \{ M(A), A \in {\cal B}_0(S)\} $ of r.v.'s indexed by sets $A \in {\cal B}_0(S) := \{ A \in  {\cal B}(S): \mu(A) < \infty\} $
of a measure space $(S, {\cal B}(S), \mu) $
with a $\sigma-$finite measure $\mu $ is called a {\it random measure with distribution $W $ and control measure $\mu$}
(where $W$ is an infinitely divisible r.v.) if for any disjoint sets $A_i \in {\cal B}_0(S)$, $i=1, \dots, n$, $n\ge 1$
r.v.'s $M(A_i), i=1,\dots,n$ are independent, $M(\cup_{i=1}^n A_i) =
\sum_{i=1}^n M(A_i)$, and
$$
\E \e^{\i \theta M(A)} = (\E \e^{\i \theta W})^{\mu(A)}, \quad \forall  A \in {\cal B}_0(S), \ \  \forall  \theta \in \R.
$$
Integration with respect to stable and general infinitely divisible random measures is discussed in \citetext{bib:ST94},
\citetext{bib:RR89} and other texts.
The process in (\ref{mixstable}) is a particularly simple case of the so-called {\it mixed stable moving averages}
introduced in \citetext{bib:SUR93}. It is stationary, ergodic, and
has $\alpha-$stable finite-dimensional distributions. The representation (\ref{mixstable}) of the limit aggregated process
holds also in the finite-variance case $\alpha =2$, yielding a stationary Gaussian process with zero mean and
covariance
$$
{\rm Cov}(\mathfrak{X}(0), \mathfrak{X}(t)) \ =\  \sigma^2 \sum_{s \le 0} \int_{[0,1)} x^{-s} x^{t-s} \Phi(\d x)\ =\
 \sigma^2 \E \Big[ \frac{a^t}{1-a^2}\Big], \quad t = 0,1, \dots,
$$
cf.\ (\ref{covX}). The mixed stable moving average in (\ref{mixstable}) can be regarded
as a limiting  ``superposition''
$\sum_{a_i \in [0,1)} {\cal X}(t; a_i)$ of independent
$\alpha-$stable AR(1) processes ${\cal X}(t; a_i) = \sum_{s \le t} a_i^{t-s} {\cal M}_{\alpha,i} (t-s) $ with $\alpha-$stable innovations
$\{{\cal M}_{\alpha,i} (s) \}$. Although each $\{{\cal X}(t; a_i)\}$ is geometrically mixing and hence short memory, the dependence
in  $\{{\cal X}(t; a_i)\}$ increases when $a_i \nearrow 1 $ approaches the ``unit root''  $a =1 $. It turns out that the limiting  ``superposition'' of these processes and the mixed moving average in (\ref{mixstable})
may have long (or short) memory, depending on the concentration of the $a_i$'s near $a =1 $, or the parameter
$\beta  $ in (\ref{mix1}).

Before  addressing the question about long memory of the infinite variance process  in (\ref{mixstable}), let us note that
the above mentioned convergence of (\ref{XN1}) to (\ref{mixstable}) does not hold in the case of negative exponent $-1 < \beta < 0$. It turns out that in the latter case, the limit aggregated process does not depend on $t$ and is an $\alpha (1+\beta)-$stable r.v.\ (random constant):
\begin{equation}\label{XN2}
\frac{1}{N^{1/(\alpha (1+\beta))}}  \sum_{i=1}^N X_{i}(t) \ \toFDD  \  \tilde Z,
\end{equation}
see Puplinskait\.e and Surgailis (\citeyear{bib:ps2010}, Prop.\ 2.3). In the case $\alpha =2$, a similar result was noted in \citetext{bib:ZA}.
Note that the normalizing exponents $1/\alpha$ in (\ref{XN1}) and $ 1/(\alpha (1+\beta))$ in (\ref{XN2})
are different and $ 1/(\alpha (1+\beta)) \to \infty $ as $\beta \to -1$.
Therefore, $\beta=0$ is a critical
point resulting in completely different limits of the aggregated process in the cases $\beta>0$ and $\beta<0$. The fact that the
limit is degenerate in the latter case can be intuitively explained as follows. It is clear that, with $\beta$ decreasing, the dependence increases
in the random-coefficient AR(1) process $\{X(t)\}$, as well as in the limiting aggregated process $\{\mathfrak{X}(t) \}$.
For negative $\beta<0$, the dependence in the aggregated process becomes
extremely strong so that the limit process is degenerate and completely dependent.
\smallskip

\noindent {\it Long memory properties of the limit aggregated process.}
Clearly, the usual definitions of long memory in terms of covariance/spectrum do not
apply in infinite variance case.  Alternative notions of long memory which are applicable to infinite-variance processes have been proposed in the literature. \citetext{bib:as83} was probably the first to rigorously study
long memory for such processes in terms of the  rate at which the bivariate characteristic function at distant
lags factorizes into the product of two univariate characteristic functions. Related
characteristics such as codifference are discussed in \citetext{bib:ST94}. Some characteristics
of dependence (covariation, $\alpha$--covariance) for stable processes expressed in terms of the spectral measure were studied in
\citetext{bib:ST94} and \citetext{bib:paul2013}. \citetext{bib:HY} defined the long-range dependence (sample
Allen variance) (LRD(SAV)) property in terms of the limit behavior of
squared studentized sample mean.

Probably, the most useful and universal definition of long/short memory was given by \citetext{bib:C84a}. Assume that
$\{Y(t), t \in \Z\}$ is a strictly stationary process series and there exist some constants $D_n \to \infty \ (n \to
\infty) $ and $B_n $ and a nontrivial stochastic  process
$J = \{J(\tau), \tau \ge 0\}\not\equiv 0$
such that
\begin{equation}
\frac{1}{D_n} \sum_{s=1}^{[n\tau]} (Y(s) - B_n) \ \toFDD  \
J(\tau). \label{distLM}
\end{equation}
According to \citetext{bib:C84a}, if the limit process $J$ has {\it independent increments}, the series  $\{ Y(t)\}$ is said to have {\it distributional short memory}, while in the converse
case when
$J$ has {\it dependent increments}, the series  $\{ Y(t)\}$ is said to have {\it distributional long memory}.

The above definition has several advantages. First of all, it does not depend on any moment assumptions since
finite and infinite variance processes are treated from the same
angle. Secondly,  according to the classical Lamperti's theorem (see \citetext{bib:L62}),
in the case of (\ref{distLM}) there exists a number $H>0$ such that the normalizing constants $A_n$ grow as
$n^H$ (modulus a slowly varying factor), while the limit random process $J$ is $H-$self-similar and has stationary
increments ($H-$sssi). A faster growth of normalization $A_n$ means stronger dependence  and therefore
$H$ is a quantitative indicator of the degree of dependence in  $\{ Y(t)\}$. The characterization of short memory
through (\ref{distLM}) is very robust and essentially reduces to L\'evy stable behavior since all sssi processes
with independent increments are stable L\'evy processes. We emphasize that the above definition of
long/short memory requires  identification of the partial sums limit which is sometimes not easy. On the other hand,
partial sums play a very important role in statistical inference, especially  under long memory, see e.g.\ \citetext{bib:book2012},
and hence finding partial sums limit is very natural for understanding the dependence structure
of a given process and subsequent applications of statistical nature.

With the above discussion in mind, the question arises what is the partial sums limit of
the aggregated process in (\ref{mixstable})? This question is answered in Puplinskait\.e and
  Surgailis (\citeyear{bib:ps2010}, Thm.\ 3.1) saying that
for $ 1< \alpha \le 2 $ and  $0< \beta < \alpha -1$,
\begin{equation} \label{LMconv}
\frac{1}{n^{1 - (\beta/\alpha)}} \sum_{t=1}^{[n\tau]} \mathfrak{X}(t) \ \ \toFDD \ \ \Lambda_{\alpha,\beta}(\tau) \
:=\ \int_{\R_+ \times \R} (\mathfrak{f}(x, \tau-s) - \mathfrak{f}(x,-s)) Z_\alpha (\d x, \d s),
\end{equation}
where
\begin{equation}\label{fdef}
\mathfrak{f}(x,t)\ := \ \begin{cases} (1- {\e}^{-xt})/x, &\text{if} \  x>0 \ \text{and} \ t>0, \\
0,&\text{otherwise},
\end{cases}
\end{equation}
and $Z_\alpha(\d x, \d s)$ is an $\alpha-$stable random measure 
on $\R_+ \times \R$ with
control measure $c_\phi x^{\beta} \d x\, \d s $ and distribution $Z$, where $Z$ is defined
at (\ref{DA}). The random process
$\Lambda_{\alpha,\beta}$ in (\ref{LMconv}) is well-defined for $ 1< \alpha \le 2$, $0<\beta<\alpha-1$ and is
$H-$sssi with self-similarity index $H = 1 - \frac{\beta}{\alpha} \in (\frac{1}{\alpha}, 1). $ Moreover,
$\Lambda_{\alpha,\beta}$ has a.s. continuous paths, $\alpha-$stable finite dimensional distributions and stationary and dependent increments. In particular,  $\Lambda_{2,\beta}$ is a fractional Brownian motion with $H = 1-\frac{\beta}{2} \in (\frac{1}{2}, 1). $
The process
$\Lambda_{\alpha,\beta}$ is also different from linear fractional L\'evy motion (see, e.g., \citetext{bib:ST94} for a detailed discussion of the  latter process).

Similarly as $\{ \mathfrak{X}(t) \}$ in (\ref{mixstable}) can be regarded
as a  ``continuous superposition'' of AR(1) processes, the limit process $\Lambda_{\alpha,\beta}$ in (\ref{LMconv})
can be regarded as a  ``continuous superposition'' of (integrated) Ornstein-Uhlenbeck processes
$z_\alpha (\tau; x)$ defined as
\begin{equation}\label{zdef1}
z_\alpha(\tau; x) \ :=\
\int_0^\tau \d u \int_{\R} \e^{-x(u-s)} \1(u>s) \d Z_\alpha(s) \ = \  \int_{\R} \big(\mathfrak{f}(x, \tau -s) - \mathfrak{f}(x, -s)\big) \d Z_\alpha  (s), \quad \tau\ge 0, \ x >0,
\end{equation}
where  $ \{Z_\alpha(s), s \ge 0\}$ is an $\alpha-$stable  L\'evy process with independent increments. Note that
for each $x > 0$, the process
$\{z_\alpha (\tau; x)\}$ is a.s.\ continuously differentiable on $\R$ and its derivative
$z'_\alpha(\tau; x) =  \d  z_\alpha (\tau; x)/\d \tau $ satisfies the Langevin equation
$$
\d z'_\alpha(\tau; x) \ =\ - x z'_\alpha(\tau; x)\d \tau  + \d Z_\alpha (\tau).
$$
In the case $\alpha=2$, $Z_2 = B $ is a usual Brownian motion and $z_2 (\tau; x) $ is a Gaussian  Ornstein-Uhlenbeck process.
The corresponding representation of $\Lambda_{2,\beta}$ in (\ref{LMconv}) may be termed  the {\it Ornstein-Uhlenbeck representation of fractional Brownian motion}, and the process $\Lambda_{\alpha,\beta},\, 1< \alpha < 2$ its stable counterpart. A related class of stationary infinitely divisible
processes with long memory is discussed in \citetext{bib:B-N2001}.

The condition $ 1< \alpha \le 2, \, 0< \beta < \alpha -1$ for the convergence in (\ref{LMconv}) is sharp and cannot be weakened.
In particular, for $\beta > \alpha -1 $ the partial sums process $n^{-1/\alpha} \sum_{t=1}^{[n\tau]} \mathfrak{X }(t) $ tends
to an $\alpha-$stable L\'evy process with independent increments (see Puplinskait\.e and
  Surgailis (\citeyear{bib:ps2010}, Thm.\ 3.1)). Therefore, large sample behaviors
of $\{\mathfrak{X }(t) \}$ for $0 < \beta < \alpha-1$ and $\beta > \alpha -1 $ are markedly different. Following the above terminology,
the aggregated AR(1) process $\{\mathfrak{X }(t) \}$ has distributional long memory if $0 < \beta < \alpha-1$ and
distributional short memory if $\beta > \alpha -1$. In the last case,  assumption   (\ref{mix1}) on the mixing
distribution can be substantially relaxed.

\citetext{bib:ps2010} also studied other characterizations of long memory of the aggregated AR(1) process $\{\mathfrak{X }(t) \}$
with infinite variance
(the LRD(SAV) property, the decay rate of codifference). A curious characterization of long memory in terms of the asymptotic behaviour
of the {\it ruin probability} in a discrete time risk insurance model with $\alpha-$stable claims $\{\mathfrak{X }(t) \}$ in (\ref{mixstable}) is obtained in \citetext{bib:KARINA2013}, following the  characterization studied in \citetext{bib:MIK00}. All
these results agree with the  above characterization in terms of the partial sums process, in the sense that $\beta = \alpha -1 $ is
the boundary between
long memory and short memory in  $\{\mathfrak{X }(t) \}$.

Finally, let us note that the general aggregation scheme of random-coefficient autoregressive processes discussed in \citetext{bib:ZA}
includes the case of common component, or {\it common innovations}. Aggregation of infinite variance AR(1) processes with common innovations
was studied in \citetext{bib:PUP09}. In this case, the limit aggregated process, say $\{\widetilde
{\mathfrak{X}}(t)\}$,  exists under normalization $1/N$ and is
a moving average with the same innovations as the original AR(1) series  (\ref{AR1}) and the moving average coefficients
given by the expectations $\E a^j = \int_{[0,1)} x^j \Phi (\d x)$. By a similar argument as in (\ref{covX}), we have
that $\E a^j \sim const\; j^{-\beta-1} $ for $\beta > -1 $. Therefore for $1/\alpha<\beta <0$, $\{\widetilde
{\mathfrak{X}}(t)\}$ is a well-defined long memory moving average with infinite variance and nonsummable coefficients
$\sum_{j=0}^\infty \E a^j = \infty $. \citetext{bib:PUP09} investigated various long memory properties of
$\{\widetilde{\mathfrak{X}}(t)\}$, including the convergence of its partial sums to a linear fractional L\'evy motion.

\section{Aggregation of  AR(1) processes: triangular array innovations}
\label{sec:aggr-ar1-proc-2}

The contemporaneous aggregation scheme of Sec.\ \ref{sec:aggr-ar1-proc} can be generalized by assuming that the
innovations depend on $N$, constituting a triangular array of i.i.d.\ r.v.'s. Such aggregation scheme was
studied in \citetext{bib:pps2012}. Let $\{ X^{(N)}_i(t) \}$, $i=1,\dots, N$ be i.i.d.\ copies of
of random-coefficient AR(1) process
\begin{equation} \label{AR2}
X^{(N)}(t)\ = \ a X^{(N)}(t-1) + \zeta^{(N)}(t), \quad t\in \Z,
\end{equation}
where $\{\zeta^{(N)} (t), \, t \in \Z \}, \, N=1,2, \dots $
is a triangular array of i.i.d.\ random variables  in the domain of  attraction of an
infinitely divisible law  $W$:
\begin{eqnarray} \label{DID}
\sum_{i=1}^N \zeta^{(N)}(i)&\toD&W,
\end{eqnarray}
and $a \in [0,1)$ is
a r.v., independent of  $\{\zeta^{(N)} (t), t \in \Z \}$.
The limit aggregated process $\{\mathfrak{X}(t), t \in \Z \} $ is defined as the limit in distribution:
\begin{equation}
\sum_{i=1}^N X^{(N)}_i(t) \ \toFDD \  \mathfrak{X}(t).\label{limaggre}
\end{equation}
Sec.\ \ref{sec:aggr-ar1-proc} corresponds to the particular case of (\ref{AR2})--(\ref{DID}), viz.,
$\zeta^{(N)}(t) = N^{-1/\alpha} \zeta(t), $ where  $ \{\zeta(t), t \in \Z \}$ are i.i.d.\ r.v.'s
in the domain of (normal) attraction of $\alpha-$stable law $W,\, 0 < \alpha \le 2. $ In particular, for $\alpha = 2$  or
$W \sim {\cal N}(0,\sigma^2)$
the last condition is equivalent to $\E \zeta(t) = 0 $ and $\sigma^2 = \E W^2 < \infty$.

One of the main results of \citetext{bib:pps2012} says that under mild additional conditions the limit in (\ref{limaggre}) exists and is written as
a mixed infinitely divisible (ID) moving average (see the terminology in \citetext{bib:RR89}):
\begin{equation}
\mathfrak{X}(t) \ = \ \sum_{s \le t} \int_{[0,1)} x^{t-s} M_{W,s}({\d}x), \qquad t \in \Z,
\label{mix}
\end{equation}
where  $\{M_{W,s}, s\in \Z\} $ are i.i.d.\ copies of an ID random measure $M_W$ on $[0,1)$ with control measure
$\Phi (\d x) = \P (a \in \d x)$ and the distribution $W$ in (\ref{DID}). Recall that the last distribution is uniquely determined
by its L\'evy characteristics $(\mu, \sigma, \pi)$ (the characteristic triplet) since
\begin{eqnarray} \label{Wchf0}
V(\theta)&:=&\log \E \e^{\i \theta W}  \ = \
 \int_{\R} (\e^{ {\i} \theta x} - 1 - {\i} \theta x \1(|x| \le 1)) \pi (\d x) - \frac{1}{2}\theta^2 \sigma^2 +  \i \theta \mu,
\end{eqnarray}
where $ \mu \in \R, \, \sigma \ge 0$ and $\pi$ is a L\'evy measure, see \citetext{bib:ST94}, \citetext{bib:SAT1999}.

\citetext{bib:pps2012} discuss distributional long/short memory properties of the aggregated process $\{\mathfrak{X}(t) \}$ in (\ref{mix})
with finite variance and a mixing density as in (\ref{mix1}).
The finite variance assumption is equivalent to  $\E (1- a)^{-1} < \infty $ and
\begin{equation} \label{Wvar}
\sigma^2_W \ :=\ {\rm Var}(W) < \infty \qquad  \Longleftrightarrow  \qquad  \int_{\R} x^2 \pi (\d x) < \infty.
\end{equation}
Note that (\ref{Wvar}) excludes the $\alpha-$stable case discussed in the previous sec. Under (\ref{Wvar})
the covariance function of
(\ref{mix}) is written as in the Gaussian case (\ref{covXG}), with $\sigma^2 $ replaced by $\sigma^2_W$. Therefore
the covariance asymptotics in (\ref{covX}) applies also for the process  $\{\mathfrak{X}(t) \}$ in (\ref{mix}), yielding
\begin{equation}\label{Svar}
{\rm Var}\Big(\sum_{t=1}^{n} \mathfrak{X}(t) \Big) \  \sim \  C n^{2-\beta}, \qquad 0< \beta < 1.
\end{equation}
From (\ref{Svar}) and the linear structure of $\{\mathfrak{X}(t) \}$
one might expect a Gaussian (fractional Brownian motion) limit behavior of the partial sums process
$S_n(\tau) = \sum_{t=1}^{[n\tau]} \mathfrak{X}(t)$.

However, as it turns out, the Gaussian scenario for $\{S_n(\tau) \}$ is valid only if $\sigma >0$ in (\ref{Wchf0}),
or the Gaussian component is present in the ID r.v.\ $W$. Else (i.e., when $\sigma =0$), the behavior of the
 L\'evy measure $\pi$ at the origin plays a dominant role. Assume that there exist $\alpha_0 >0$ and
$c^\pm_0 \ge 0, c^+_0 + c^-_0 >0$ such that
\begin{equation} \label{limPi}
\lim_{x \searrow 0} x^{\alpha_0} \pi(x, \infty) = c^+_0, \qquad  \lim_{x \searrow 0} x^{\alpha_0 } \pi(-\infty, -x) = c^-_0.
\end{equation}

It is proved in \citetext{bib:pps2012}
that under conditions (\ref{mix1}), (\ref{Wvar}), and (\ref{limPi}),  partial
sums $\{S_n(\tau) \}$ of $\{\mathfrak{X}(t) \}$ in (\ref{mix}) may exhibit at least four
different limit behaviors, depending on parameters
$\beta, \sigma, $ and $\alpha_0$.
The four parameter regions and the limit behaviors in the f.d.d.\ sense
are described in (i)--(iv) below.
\begin{description}
\item [{\rm (i)}]  \ $0< \beta < 1, \  \sigma >0. $  In this region,  $\ n^{(\beta/2) - 1}S_n(\tau)\ $ tends to a fractional Brownian motion with Hurst parameter $H = 1 - (\beta/2).  $
\item [{\rm (ii)}] \ $0 < \beta < 1, \ \sigma = 0, \ 1+ \beta < \alpha_0 < 2. $  In this region,  $\ n^{(\beta/\alpha_0) - 1}S_n(\tau)\ $ tends to the
$\alpha_0-$stable self-similar process $\Lambda_{\alpha_0, \beta}$ defined in (\ref{LMconv}).
\item [{\rm (iii)}] \ $0 < \beta < 1, \ \sigma = 0, \ 0 < \alpha_0 <  1 + \beta. $
In this region,  $\ n^{-1/(1+\beta)}S_n(\tau)\ $ tends to a
$(1+\beta)-$stable   L\'evy process with independent increments.
\item [{\rm (iv)}] \ $\beta > 1.$
In this region,  $\ n^{-1/2}S_n(\tau)\ $ tends to
a Brownian motion.
\end{description}

See \citetext{bib:pps2012} for  precise formulations.
Accordingly, the process  $\{\mathfrak{X}(t) \} $ in (\ref{mix}) has distributional long memory in regions (i) and (ii) and distributional
short memory in regions (iii) and (iv).
As $\alpha_0 $ increases from $0$ to $2$, the
L\'evy measure in (\ref{limPi}) increases its ``mass'' near the origin, the limiting case $\alpha_0 = 2$
corresponding to $\sigma >0$ or a positive ``mass''
at $0$. We see from (i)--(ii) that distributional long memory is related to $\alpha_0 $ being large enough, or
small jumps of the random
measure $M_W$ having sufficient high intensity.  Note that  the critical exponent $\alpha_0 = 1 + \beta$ separating the
long and short memory ``regimes'' in (ii) and (iii)
decreases  with $\beta, $ which is quite natural since  smaller $\beta $ means the mixing distribution putting more weight
near the unit root $a=1$.

Let us note that an $\alpha-$stable limit behavior of partial sums of stationary finite variance processes is not unusual
under long memory. See, e.g., \citetext{WTLW:1997}, \citetext{bib:MRRS02}, \citetext{bib:LS2003} and the references herein. On the other hand,
these papers focus on heavy-tailed duration models in which case
the limit $\alpha-$stable process has {\it independent increments} as a rule. The situation when an infinite variance limit process
with {\it dependent increments} arises from partial sums of a finite variance process as in (ii) above seems rather new.

\section{Joint temporal and contemporaneous aggregation of  AR(1) processes }\label{sec:joint}

The aggregation procedures discussed in the previous sec.\ extend in a natural way to the (large scale)
joint temporal and contemporaneous aggregation. In the latter frame, we are interested in  the limit behavior of the double sums
\begin{equation} \label{aggreSum}
S_{N,n}(\tau) \ := \  \sum_{i=1}^N \sum_{t=1}^{[n\tau]} X_i(t), \quad \tau \ge 0,
\end{equation}
where
$\{X_i(t)\}, i=1, \dots, N$ are the same random-coefficient AR(1) processes as in (\ref{XN}). The sum in  (\ref{aggreSum})
represents joint temporal and contemporaneous aggregate of $N$ individual AR(1) evolutions
at time scale $n$. The main question is the joint aggregation limit $\lim_{N, n \to \infty} A^{-1}_{N,n} S_{N,n}(\tau)$,
in distribution, where $A_{N,n}$ are some normalizing constants and both
$N$ and $n$ increase to infinity,  possibly at different rate. This question was studied in \citetext{pil2013}.
The last paper also discussed the
iterated limits of $A^{-1}_{N,n} S_{N,n}(\tau) $
when first $n \to \infty $ and then
$N \to \infty $, or vice-versa.  Remark that the discussion in the previous sec.\ refers to the latter iterated limit as   $N\to \infty$ first,
followed by
$n\to \infty$. Similar questions
for some network traffic models were studied in \citetext{WTLW:1997}, \citetext{bib:MRRS02}, \citetext{gaig2003}, \citetext{pip2004},
\citetext{domb2011} and other papers. In these papers, the role of AR(1) processes $\{X_i(t)\}$ in (\ref{aggreSum}) play
independent ON/OFF processes or M/G/$\infty$ queues with heavy-tailed  activity periods.

Let us describe the main results in \citetext{pil2013}. They refer to the random-coefficient AR(1) process with i.i.d.\ innovations having
zero mean and variance $\sigma^2 < \infty $, and a mixing density as (\ref{mix1}).
Let $N, n $ increase simultaneously so as
\begin{equation}\label{Nninc}
\frac{N^{1/(1+\beta)}}{n} \ \to \  \mu \, \in \, [0,\infty],
\end{equation}
leading to the following three cases: 
\begin{equation}\label{cases}
\text{Case (j)}: \ \mu = \infty, \qquad \text{Case (jj)}: \  \mu = 0, \qquad  \text{Case (jjj)}: \ 0< \mu < \infty.
\end{equation}
Following the terminology in \citetext{bib:MRRS02} and \citetext{pip2004}, we call Cases (j), (jj), and (jjj)
the ``fast growth condition'' , the ``slow growth condition'',  and the ``intermediate growth condition'', respectively,
since they reflect how fast $N$ grows with $n$.  The main result of \citetext{pil2013}
says that under (\ref{Nninc}), the ``simultaneous limit'' of
$S_{N,n}(\tau)$ exist and are different in all three Cases (j)--(jjj).

\smallskip

\noindent {\it Case (j)} (the ``fast growth condition''): For any $0< \beta < 1$,
\begin{equation}
N^{-1/2} n^{-1 + (\beta/2)} S_{N,n}(\tau)\ \toFDD\  B_{1-\beta/2}(\tau), \label{casej}
\end{equation}
where $\{B_{1-\beta/2}(\tau)\}$ is a fractional Brownian motion with $H = 1-\beta/2 \in (1/2, 1)$.

\smallskip

\noindent {\it Case (jj)} (the ``slow growth condition''): For any $-1< \beta < 1$,
\begin{equation}
N^{-1/(1+\beta)} n^{-1/2} S_{N,n}(\tau)\ \toFDD\  {\cal W}_{\beta}(\tau), \label{casejj}
\end{equation}
where $\{{\cal W}_{\beta}(\tau)\}$ is a sub-Gaussian $(1+\beta)-$stable process defined as
${\cal W}_{\beta}(\tau) = W_\beta^{1/2} B(\tau), \, \tau \ge 0$, where $W_\beta >0$ is a $(1+\beta)/2-$stable totally skewed
r.v. and $\{B(\tau), \tau \ge 0\}$ is a standard Brownian motion independent of $W_\beta$ (see, e.g., \citetext{bib:ST94}).

\smallskip

\noindent {\it Case (jjj)} (the ``intermediate growth condition''): For any $-1< \beta < 1$
\begin{equation}
N^{-1/(1+ \beta)} n^{-1/2} S_{N,n}(\tau)\ \toFDD\  \mu^{1/2}
{\cal Z}_{\beta}(\tau/\mu),  \label{casejjj}
\end{equation}
where the limit process $\{{\cal Z}_{\beta} (\tau), \tau \ge 0 \} $ is defined through finite-dimensional characteristic function:
\begin{eqnarray}
\E \exp \Big\{ \i \sum_{j=1}^m \theta_j {\cal Z}_{\beta} (\tau_j)  \Big\} \label{calZfd}
&=&\exp \Big\{ c_\phi \int_{\R_+}  \big(\e^{- (\sigma^2/2) \int_{\R}  \big(\sum_{j=1}^m \theta_j (\mathfrak{f}(x, \tau_j -s) - \mathfrak{f}(x, -s)) \big)^2 \d s} -1 \big) x^\beta \d x  \Big\},
\end{eqnarray}
where $\theta_j \in \R, \tau_j \in \R_+, j=1, \dots, m, \, m\ge 1,$ and $\mathfrak{f}$ is given in (\ref{fdef}).

\smallskip

Let us give some comments on the above results. The convergence in  Case (j) (\ref{casej})  is very natural in view of
(\ref{LMconv}) (with $\alpha =2$ and $\Lambda_{2,\beta} = B_{1-\beta/2}$ a fractional Brownian motion). Indeed, under
the ``fast growth condition'' we expect that $N^{-1/2} n^{-1 +
  (\beta/2)} S_{N,n}(\tau)$ can be approximated by $ n^{-1 + (\beta/2)}\sum_{t=1}^{[n\tau]} \mathfrak{X}(t)$, which converges to $B_{1-\beta/2}(\tau)$ as it happens in the case
when $N \to \infty$ followed by $n \to \infty$. The limit in Case (jj) (\ref{casejj}) can be also easily explained since
``individual'' partial sums  $ \sum_{t=1}^{[n\tau]} X_i(t)$ behave as Brownian motions with random variances $(1-a_i)^{-2}$ having infinite
expectation and a heavy tailed distribution with tail parameter $(1+ \beta)/2 \in (0,2)$:  $\P ((1-a_i)^{-2} > x) = \P (a_i > 1 - 1/\sqrt{x})
\sim C x^{- (\beta +1)/2}, \, x \to \infty.  $ The sum of such independent ``random-variance'' Brownian motions
behaves as a sub-Gaussian process  ${\cal W}_{\beta}$ in  (\ref{casejj}). Particularly interesting is the limit process ${\cal Z}_\beta$ arising
under ``intermediate scaling'' in Case (jjj). It is shown in \citetext{pil2013} that ${\cal Z}_\beta$ admits  a  stochastic integral
with respect to a Poisson random measure  
on the product space $\R \times C(\R)$ with mean
$\psi_1 x^\beta \d x \times \P_B$, where $\P_B$ is the Wiener measure on $ C(\R)$ and
enjoys several
``intermediate'' properties between the limits in (j) and (jj).  According to (\ref{calZfd}), ${\cal Z}_\beta$ has infinitely divisible
finite-dimensional distributions and stationary increments, but is neither self-similar nor stable. For $0<\beta < 1 $  the process ${\cal Z}_\beta$
has finite variance
and the covariance equal to that of a fractional Brownian motion.
These results can be compared to \citetext{bib:MRRS02}, \citetext{gaig2003}, \citetext{pip2004}, \citetext{domb2011}.
In particular, \citetext{bib:MRRS02} discuss
the ``total accumulated input''  ${\cal A}_{n,N}(\tau) := \int_0^{n \tau} \sum_{i=1}^N W_i (t) \d t $ from $N$ independent ``sources''
$\{ W_i(t), i=1, \dots, N; t \ge 0 \}$
at time scale $n$. The aggregated inputs are i.i.d. copies of ON/OFF process
$\{ W_t, t \ge 0 \}, $ alternating between 1 and 0 and taking value 1 if $t$ is in an ON-period and 0 if $t$
is in an OFF-period, the ON- and OFF-periods forming a stationary renewal process having heavy-tailed
lengths with respective tail parameters $\alpha_{\rm on}, \alpha_{\rm off} \in (1,2), \, 
\alpha_{\rm on} < \alpha_{\rm off}$. The role of condition (\ref{Nninc}) is played in the above paper by
\begin{equation}\label{Nninc1}
\frac{n}{N^{\alpha_{\rm on} -1}} \ \to \  \mu \, \in \, [0,\infty],
\end{equation}
leading to the three cases analogous to (\ref{cases}):
\begin{equation}\label{cases1}
\text{Case (j')}: \ \mu = 0, \qquad \text{Case (jj')}: \  \mu = \infty, \qquad  \text{Case (jjj')}: \ 0< \mu < \infty.
\end{equation}
The limit of (normalized) ``input''  ${\cal A}_{n,N}(\tau)$ in Case (j') (the ``slow growth condition'' )
and Case (jj') (the ``fast growth condition'' ) was obtained in \citetext{bib:MRRS02}, as an $\alpha-$stable L\'evy process and
a fractional Brownian motion, respectively. The ``intermediate'' limit in Case (jjj') was identified in
\citetext{gaig2003},  \citetext{gaig2006}, \citetext{domb2011}
who showed that this process can be regarded as a ``bridge'' between
the limiting processes in Cases (j') and (jj'), and
can be represented as a stochastic integral with respect to a Poisson random measure on $\R_+ \times \R $.
A common feature to the above research and  \citetext{pil2013} is the fact the
partial sums of the individual processes with finite variance tend to an infinite variance process,
thus exhibiting an increase of variability. See also \citetext{bib:LS2003}, \citetext{lei2005}.
These analogies
raise interesting open questions about extension of the joint temporal-contemporaneous aggregation scheme to general independent
processes with covariance long memory and stable behavior of partial sums.

\section{Aggregation of autoregressive random fields}\label{sec:fields}

The idea of aggregation naturally extends to spatial autoregressive models (\citetext{bib:lav2007},
\citetext{bib:lav2011}, \citetext{pup2013}, \citetext{bib:lls2013}). Following \citetext{pup2013}, consider a nearest-neighbor autoregressive random field $\{X(t,s)\}$  on $\Z^2$ satisfying
the difference equation
\begin{eqnarray} \label{modelNN}
X(t,s)&=&\sum_{|u|=|v|=1} a(u,v) X(t+u,s+v) + \zeta(t,s), \quad (t,s) \in \Z^2,
\end{eqnarray}
where 
$\{\zeta(t,s), (t,s) \in \Z^2 \} $ are i.i.d. r.v.'s whose generic distribution $\zeta  \in D(\alpha)$ belongs
to the domain of (normal) attraction of $\alpha-$stable law,
$1< \alpha \le  2$, and
$a(t,s) \ge 0, \, |t|=|s|=1$ are {\it random} coefficients, independent of $\{\zeta(t,s) \}$ and satisfying the  condition
$A := \sum_{|t|=|s|=1} a(t,s)<1$ a.s. for the existence of a stationary solution of (\ref{modelNN}).
The stationary solution of (\ref{modelNN}) is given by the convergent series
\begin{eqnarray} \label{stationary0}
X(t,s)&=&\sum_{(u,v) \in \Z^2 } g(t-u,s-v, a)\zeta(u,v), \qquad (t,s) \in \Z^2,
\end{eqnarray}
where $g(t,s,a), \, (t,s) \in \Z^2, a = (a(t,s), |t|=|s| =1) $, is the (random) lattice Green function solving the equation $g(t,s,a)-  \sum_{|u|=|v|=1}a(u,v) g(t+u,s+v, a) = \delta(t,s), $
where $\delta (t,s)$ is the delta function. Under the condition $\E (1-A)^{-1} < \infty $,
the series in (\ref{stationary0}) converges unconditionally in $L^p $ for any $1< p < \alpha $ (\citetext{bib:lls2013}).
In the finite variance case $\alpha = 2$,
the stationary solution (\ref{stationary0}) can be defined via spectral representation:
\begin{eqnarray} \label{stationary00}
X(t,s)&=&\int_{[-\pi, \pi]^2} \e^{\i (tx + sy)}\hat g(x,y, a) Z(\d x, \d y),
\end{eqnarray}
where $\hat g(x,y, a) = \big(1- \sum_{|t|=|s|=1} a(t,s) \e^{\i (xt + ys)}\big)^{-1} $ is the Fourier transform of $g(t,s,a)$
and $Z(\d x, \d y)$ is the random meausure sarisfying $\int_{[-\pi, \pi]^2}\e^{\i (tx + sy)} Z(\d x, \d y) = \zeta(t,s).$
The spectral density of (\ref{stationary00}) is written similarly to one-dimensional case
\begin{equation}\label{sdensity}
f(x,y) = \frac{\sigma^2}{(2\pi)^2}\, \E |\hat g(x,y, a)|^2, \quad (x,y) \in [-\pi, \pi]^2, \quad \sigma^2 = \E \zeta^2.
\end{equation}
Let $\{X_i(t,s)\}, \, i=1,2, \ldots $ be independent copies of (\ref{stationary0}). The aggregated field
$\{ \mathfrak{X}(t,s)\}$ is defined as the limit in distribution:
\begin{eqnarray} \label{aggreNN}
N^{-1/\alpha} \sum_{i=1}^N X_i(t,s)&\limfdd&\mathfrak{X}(t,s), \quad (t,s) \in \Z^2.
\end{eqnarray}
Under mild additional conditions, \citetext{pup2013} prove that
the limit  in (\ref{aggreNN}) 
exists and is written as a stochastic integral
\begin{eqnarray} \label{aggremix0}
\mathfrak{X}(t,s)
&=&\sum_{(u,v) \in \Z^2} \int_{{\mbf A}} g(t-u,s-v, a) M_{u,v}(\d a), \qquad (t,s) \in \Z^2,
\end{eqnarray}
where $\{ M_{u,v}(\d a), \, (u,v) \in \Z^2 \} $ are i.i.d. copies of an
$\alpha-$stable random measure $M$  on ${\mbf A} := \{ a_{t,s} \in [0,1), \,
\sum_{|t|=|s|=1}a_{t,s} < 1 \} \subset  \R^4 $ with control measure equal to the (mixing) distribution  $\Phi$  of the random vector $a = (a(t,s), |t|=|s|=1) $ taking values in ${\mbf A}$.
The random field $\{ \mathfrak{X}(t,s)\}$
in (\ref{aggremix0}) is $\alpha-$stable and
a particular case of mixed stable moving-average fields introduced in \citetext{bib:SUR93}.

It is not surprising that dependence properties of the random field $\{ \mathfrak{X}(t,s)\}$ in (\ref{aggremix0})
strongly depend on the concentration of $\Phi$ near
the ``unit root boundary'' $A = \sum_{|t|=|s|=1} a_{t,s} =1$ but also on the form of the autoregressive operator in
(\ref{modelNN}). \citetext{pup2013} assume
that the `angular coefficients'  $0\le p(t,s) := a(t,s)/A,  \, \sum_{|t|=|s|=1} p(t,s) =1 $
are nonrandom and discuss   the following three equations
\begin{eqnarray}
X(t,s)&=&\frac{A}{2}\big(X(t- 1,s) + X(t,s- 1)\big)  +
\zeta(t,s),
\label{2N} \\
X(t,s)
&=&\frac{A}{3} \big(X(t-1,s) +  X(t,s+1)  + X(t,s-1)\big) +  \zeta(t,s), \label{3N}
\\
X(t,s)&=&\frac{A}{4}\big(X(t-1,s) +  X(t+1,s)  + X(t,s+1) + X(t,s-1)\big) + \zeta(t,s) \label{4N}
\end{eqnarray}
termed the
2N, 3N and 4N models, respectively (N standing for "Neighbour"), with a random
`radial coefficient'  $A \in [0,1)$ having a regularly varying probability density $\phi$  at $a=1$:
\begin{eqnarray} \label{mixdensity}
\phi(a)&\sim&\phi_1 (1-a)^\beta, \quad \qquad a \nearrow  1, \quad\exists \, \phi_1 >0, \
0< \beta < \alpha -1, \, 1 < \alpha \le 2.
\end{eqnarray}
Stationary solution of the above  equations in all three cases is given by (\ref{stationary0}), the Green function
being written as
\begin{eqnarray}\label{green00}
g(t,s,a)&=&\sum_{k=0}^\infty a^{k} p_k(t,s), \qquad (t,s) \in \Z^2, \qquad a \in [0,1),
\end{eqnarray}
where $p_k(t,s) = \P(W_k = (t,s) |W_0 = (0,0))$ is the $k-$step probability of the nearest-neighbor random walk $\{ W_k, k =0,1,\cdots \} $
on the lattice $\Z^2$ with one-step transition
probabilities given by
$p(t,s) = 1/2 \,  (t + s =1, t\ge 0 , s \ge 0)  $ (2N), $p(t,s) = 1/3 \,  (t + |s| =1, t\ge 0)  $ (3N) and
$p(t,s) = 1/4 \,  (|t| + |s| =1)  $ (4N), respectively.

Studying long memory properties of mixed moving average random fields $\{\mathfrak{X}_i(t,s)\}, i=2,3,4$  in (\ref{aggremix0}) corresponding
to model equations  (\ref{2N})--(\ref{4N})
requires  the control of
the Green functions in (\ref{green00}) as $|t|+|s| \to \infty $ and $a \to 1$ {\it simultaneously} since for any $a_0 < 1 $ fixed
the $g(t,s,a)$'s, $a < a_0$ decay exponentially fast with  $|t|+|s| \to \infty $. Moreover, the 2N and 3N models exhibit strong
anisotropy, characterized by a markedly different scaling behavior from the 4N model. In the Gaussian case $\alpha =2$ the random
fields (\ref{aggremix0})  are completely determined by their spectral
density $f(x,y)$ in (\ref{sdensity}) and the scaling properties of (\ref{aggremix0})  essentially
reduce to the low frequency asymptotics of $f(x,y)$ as $(x,y) \to (0,0)$. For the lattice isotropic 4N model and $\alpha =2$
the asymptotics of $f(x,y) \equiv f_4(x,y)$ under (\ref{mixdensity}) was studied by \citetext{bib:lav2011} (see also \citetext{bib:azo2009}):
$$
f_4(x,y)\ = \ const \int_0^1 \frac{1}{ |1- (a/4)\sum_{|t|=|s|=1} \e^{\i (xt + ys)}\big|^2} \phi(a) \d a   \,
\sim \, \frac{const}{(x^2 + y^2)^{1-\beta}},  \quad (x,y) \to (0,0).
$$
This shows that the large-scale limit of the Gaussian field $\{\mathfrak{X}_4(t,s)\}$ is a fully isotropic self-similar
generalized random field on $\R^2$.  However, for the 2N model the behavior of the spectral density is different
$$
f_2(x,y)\ = \ const \int_0^1 \frac{1}{ |1- (a/2)(\e^{\i xt} + \e^{\i  ys})\big|^2} \phi(a) \d a   \,
\sim \, \frac{const}{|x+y|^{1-\beta}}
K\Big(\frac{(x-y)^2}{|x+y|}\Big),  \quad (x,y) \to (0,0),
$$
where $K (u), u \ge 0 $ is a bounded continuous function on $[0, \infty)$ with $K(u) \sim const/u^{1-\beta}, u \to \infty $
(\citetext{bib:lls2013}). The above relations can be rewritten as
$$
\lim_{\lambda \to 0}\lambda f_4(\lambda^{1/H} x,\lambda^{1/H} y) = h_4(x,y) := \frac{const}{(x^2 + y^2)^{1-\beta}}, \quad H = 2(1-\beta)
$$
and
$$
\lim_{\lambda \to 0}\lambda f_2(\lambda^{1/H_1}b_1 x + \lambda^{1/H_2}b_2 y,\lambda^{1/H_1}b_3 x + \lambda^{1/H_2}b_4 y ) = h_2(x,y) :=
\frac{const}{|x|^{1-\beta}} K\Big(\frac{2y^2}{|x|}\Big), \quad H_1 = 1-\beta, \quad H_2 = 2(1-\beta),
$$
where $B =  \left(\begin{array}{cc}
b_1 & b_2 \\
b_3 & b_4 \\
\end{array}
\right) =    \left(\begin{array}{cc}
1 & -1 \\
1 & 1 \\
\end{array}
\right)  $ is non-degenerated $2\times 2-$matrix. Note that the limit function $h_2$ is non-degenerated, in the sense that it depends on
both coordinates $x $ and $y$, and that the scaling limit for $f_2$ involves two different scaling exponents $H_1 \neq H_2$, as well
as a linear transformation of $\R^2$ with non-degenerated matrix $B$. \citetext{bib:lls2013} argue that the above behavior of $f_2$ is characteristic
to anisotropic long memory, in contrast to the isotropic long memory behavior of $f_4$.

\citetext{pup2013} define {\it anisotropic distributional long memory} through scaling behavior, or partial sums limits
\begin{equation} \label{Xsum}
n^{-H_1} \sum_{(t,s) \in \tilde K}\mathfrak{X}(t,s) \ \limfdd \ V(x,y), \qquad (x,y) \in \R^2_+ = \{(x,y) \in \R^2: x,y >0\}
\end{equation}
on {\it incommensurate} rectangles
$\tilde K = K_{[nx, n^{H_1/H_2}y]} := \{ (t,s)  \in \Z^2: 1\le t \le nx, 1\le s \le n^{H_1/H_2}y \} $ with sides growing at different
rates  $O(n)$ and $O(n^{H_1/H_2}), \, H_1 \neq H_2$. 

\begin{figure}[h]
  \centering
\includegraphics[width=.4\textwidth,height=.3\textheight]{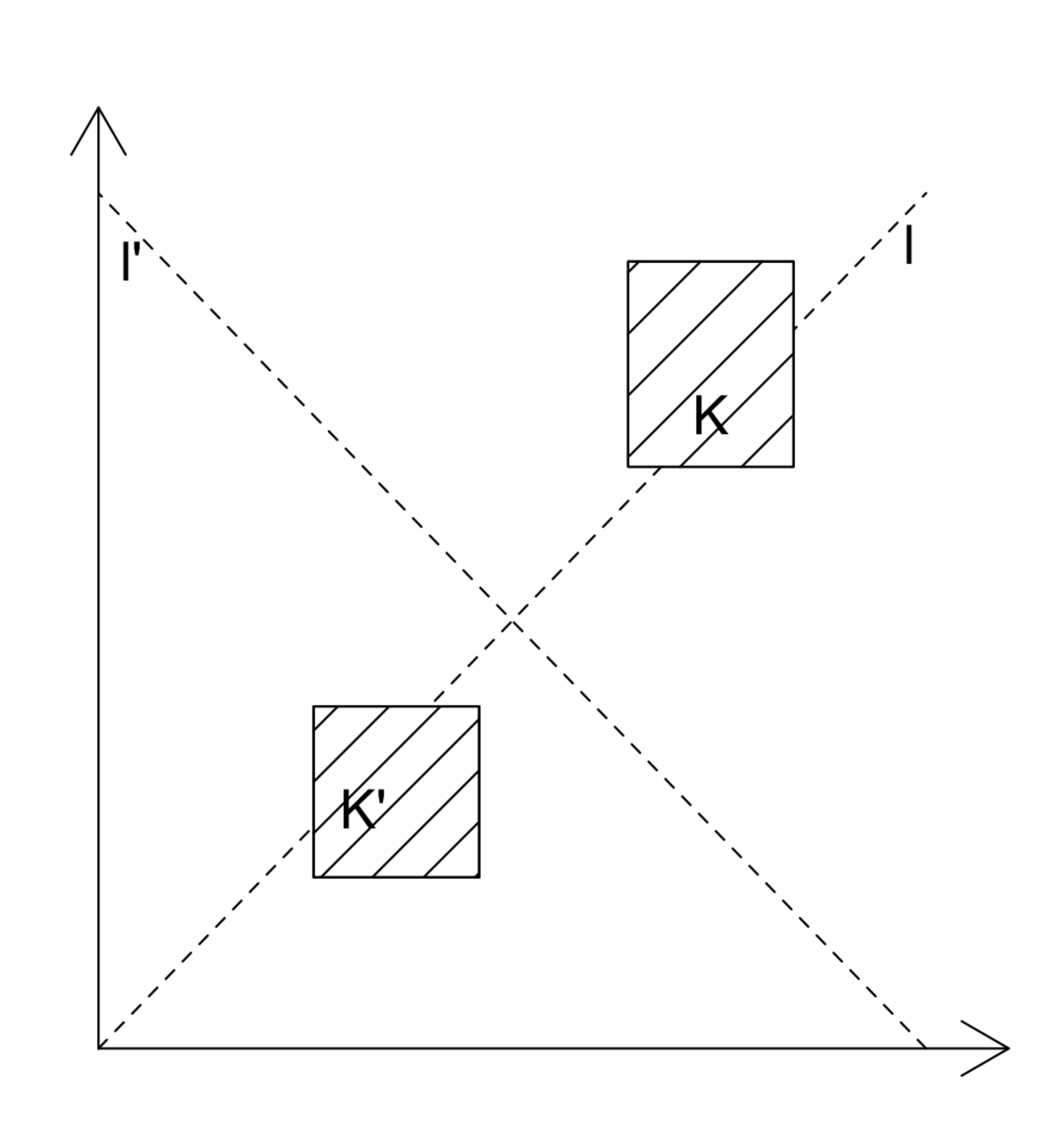}
  \caption{Location of rectangles $K$ and $K'$.}
\label{fig:I}
\end{figure}

The limit random field $\{V(x,y)\}$ in (\ref{Xsum}) is assumed to have
{\it dependent increments}, in the following sense. Given a continuous-time random field  $V = \{V(x,y)\}$, the {\it increment of $V$ on
rectangle} $K = K_{(u,v), (x,y)} = \{ (z,w) \in \R^2: u < z \le x, u < w \le y \}$ (with sides parallel to the coordinate axes) is
defined as the (double) difference:
$$
V(K) := V(x,y) - V(u,y) - V(x,v) + V(u,v).
$$
If $\ell = \{ (x,y)\in R^2:  ax + by = c \} \subset \R^2$ is a line, we say that {\it rectangles $K $ and $K'$ are separated by $\ell $} if they
lie on different sides of $\ell$. Given a line $\ell = \{ (x,y)\in \R^2:  ax + by = 0 \}$, we say that
{\it $V$ has independent increments in direction $\ell$} if for any orthogonal line $\ell' = \{ (x,y)\in \R^2:  a'x + b' y = c' \},
\ell'  \bot \ell $ and any two rectangles $K, K' $ separated by $\ell'$, increments $V(K)$ and $V(K')$ are independent r.v.'s;  else
{\it $V$ has dependent increments in direction $\ell$} (see Fig.\ \ref{fig:I}). Finally, we say that a {\it random field $V$ has dependent increments} if $V$
has dependent increments in each direction.  The above definition of anisotropic distributional long memory is contrasted in \citetext{pup2013} with
that of {\it isotropic distributional long memory}, with the only difference that in the latter case,
(\ref{Xsum}) is supposed to hold with $H_1 = H_2$, in other words, the rectangles
$K_{[nx, ny]}$ grow proportionally in each direction. In both cases, the limit random field satisfies the scaling
relation
\begin{equation}\label{OS1}
\{ \lambda V(x,y) \}   \ \eqfdd \     \{ V(\lambda^{1/H_1} x, \lambda^{1/H_2} y) \},
 \qquad \forall \, \lambda > 0.
\end{equation}
(\ref{OS1}) is a particular case of {\it operator scaling random field} (osrf) introduced in \citetext{bib:bier2007}; for $H_1 = H_2 $ (\ref{OS1})
agrees with {\it self-similar random field} (ssrf). \citetext{pup2013} show that the aggregated 2N and 3N $\alpha-$stable random
fields $\{\mathfrak{X}_2(t,s)\}$ and $\{\mathfrak{X}_3(t,s)\}$ satisfy the above definition of anisotropic long memory with
$H_1 = \frac{(1/2) + \alpha  -\beta}{\alpha}, H_2 = 2H_1$ for any $1< \alpha \le 2, 0< \beta < \alpha -1$, and identify
the limiting osrf's $V_2 $ and $V_3$  as stochastic integrals with respect to an $\alpha-$stable random measure with integrands
involving the limiting Green functions $h_2$ and $h_3$. On the other hand, the 4N field $\{\mathfrak{X}_4(t,s)\}$ is proved
to satisfy the isotropic long memory property with $H_1=H_2 = 2(\alpha - \beta)/\alpha $ and a limiting field $V_4$ involving
the limit Green function $h_4(t,s,z) = \lim_{\lambda \to \infty} g_4([\lambda t], [\lambda s], 1 - z/\lambda^2)$. The Green functions $h_i, i=2,3,4$ have a classical form of potentials of one-dimensional
heat equation and the Helmholtz equation in $\R^2$ (see \citetext{pup2013}).

\section{Disaggregation}\label{sec:disaggr}

Aggregation by itself inflicts a considerable loss of information about the evolution of individual ``micro-agents'', the latter being largely
determined by the mixing distribution.
The 
{\it  disaggregation problem} is to recover the  ``lost information'', or the mixing distribution from the spectral density or other characteristics of the aggregated series. As a
first step in this direction  we need to identify the class of spectral densities which arise from aggregation of short memory (SM) processes
with random parameters. \citetext{bib:DCF06}, \citetext{bib:DCO01} obtained an analytic characterization of the class of spectral densities
which can be written as mixtures of  infinitely differentiable SM  spectral densities. Further results in this direction were obtained
in \citetext{bib:CE}.
For example, the well-known FARIMA(0,$d$,0) spectral  density $f(y) =  (2\pi)^{-1} |2\sin (|y|/2)|^{-2d}$
can be written as the mixture  $f(y) = \int_0^1 \phi(x) |1 - x \e^{\i y}|^{-2} \phi(x) \d x $  of the AR(1) processes
corresponding to  the mixing density
$$\phi(x)\ =\ C(d) x^{d-1} (1-x)^{1-2d} (1+x). 
$$
See \citetext{bib:CE} for this  and some other examples.

The above disaggregation problem naturally leads  to  statistical problems such as estimation of the mixing distribution
from observed data. The observations  may come either from the (limit) aggregated process, or
from observed individual processes (sometimes also called panel data). Below
we review several methods of estimation of $\phi$ in the autoregressive
aggregation scheme described in  Sec.\ \ref{sec:introduction}.
Note that \eqref{covXG} means that the covariances of the
limit aggregated process  are nothing
else but the moments of the finite measure having density
$(1-x^2)^{-1}\phi (x)$. 
In fact, being supported by $[0,1)$,
this measure is uniquely determined by its moments, so
that the statistical problem of recovering the mixture density from observations $\mathfrak{X}(1), \dots, \mathfrak{X}(n)$ with finite variance
is relevant. The last problem is much harder and is still open
for the infinite variance aggregated process $\{\mathfrak{X}(t)\}$
defined in (\ref{mixstable}) (Sec.\  \ref{sec:aggr-ar1-proc}).

\smallskip

\noindent {\it Estimation from panel data}.  Consider a panel of $N$
independent AR(1) processes, each of length $n$. \citetext{bib:ROB78}
and \citetext{bib:BSG10} give estimates of
$\phi$  under the assumption that  $\phi$ belongs to some parametric family. The
main example is the family of Beta-type distribution of the form
\begin{equation}\label{dp_2}
\phi_{p,q} (x)\ =\ \frac{2}{B(p,q)}\, x^{2p-1}(1-x^2)^{q-1}, \quad x \in [0,1),\ p>1, q>1.
\end{equation}
 In this context the estimation of $\phi$ is reduced to the estimation
 of parameter $(p,q)$. We outline briefly the two approaches.

\citetext{bib:ROB78} suggested to use the classical method of
moments using the following relation between the moments of
$\phi$  and  the auto-covariances $\gamma(k) =\E X_j(0)X_j(k)$ \eqref{covX} of individual AR(1) processes,
\begin{equation}
\mu_k\ =\ \int_0^1 x^k \phi(x) \d x\ =\
\frac{\gamma(k)-\gamma(k+2)}{\gamma(0)-\gamma(2)}.\label{cle}
\end{equation}
From panel data, the $\gamma(k)$'s can be estimated
as $$\frac{1} {(n-k+1)N} \sum_{t=1}^{n-k}\sum_{j=1}^N
X_j(t)X_{j}(t+k).$$  \citetext{bib:ROB78}  proved the asymptotic normality of the corresponding estimates of moments $\mu_k, 0\le k  \le n$
when $N$ goes to infinity and $n$ is fixed.

\citetext{bib:BSG10} proposed an alternative method based on
the maximum likelihood estimate (MLE) calculated from estimated
observations. The  unobserved coefficients $a_i$ of autoregressive processes $\{X_i(t)\}$, $i=1,\dots,
N$, are estimated  by a truncated version of lag-one autocorrelation
\begin{equation*}
\widehat{a}_{i,n,h}\ = \ \min\{\max\{\widehat{a}_{i,n},h\},1-h\} \ \ {\text {with} }\ \
\widehat{a}_{i,n} \ =\
\frac{\sum_{t=1}^n X_i(t)X_i(t-1)}{\sum_{t=1}^nX_i^2(t)} , \quad h> 0.
\end{equation*}
In this way we obtain $N$ "pseudo" observations $\widehat{a}_{1,n,h}$,
$\widehat{a}_{2,n,h}$, \dots, $\widehat{a}_{N,n,h}$ of r.v.\
$a$. Then,  the parameters $p$ and $q$ in (\ref{dp_2}) are estimated by maximizing the
likelihood, viz.\ $(\hat p,\hat q)=\arg \max_{p,q} \prod_{i=1}^N \phi_{p,q}
(\widehat{a}_{i,n,h})$.  \citetext{bib:BSG10} proved the convergence
in probability  of the above MLE estimate and its
asymptotic normality with the convergence rate $\sqrt{N}$  under the
following conditions on the sample sizes and the truncation parameter $h$:
$n,N\to \infty $, $h\to 0$,
$(\log h)^2/\sqrt{ N} \to 0$,    $\sqrt{N} h^{\min(p,q)} \to 0$ and
$\sqrt{N} h^{-2} n^{-1}\to 0$.
\smallskip

\noindent {\it Estimation from the limit aggregated process.}
Assume that  a sample $\mathfrak{X}(1), \dots, \mathfrak{X}(n)$    of size $n$ is observed from
the limited aggregated process.
Under a parametric assumption about $\phi$, the  estimator proposed by
\citetext{bib:ROB78} can be easily adapted to this context, because
the limit aggregated process and the individual AR(1) have  the same covariances.

\citetext{bib:LOPV06}, \citetext{bib:CLP10}  use   the relation \eqref{cle} to construct
a non-parametric estimate of $\phi$  under the assumption that
\begin{equation}\label{LOPVdensity}
\phi(x) \ =\ (1-x)^{\beta_1}(1+x)^{\beta_2}\psi(x), \qquad   \beta_1>0, \ \beta_2>0,
\end{equation}
where $\psi(a)$ is continuous on $[-1, 1]$ and does not vanishes at
$+1$, $-1$, implying  $\E(1-a^2)^{-1} < \infty$. The above-mentioned nonparametric estimator is based on
the expansion of the mixing density  in the orthonormal basis of Gegenbauer's
polynomials in the space 
$L^2(w^{(\alpha)})$ with weight function $w^{(\alpha)}(x)=(1-x^2)^\alpha$ and $\alpha
>-1$.  The estimate is defined as
\begin{equation}\label{dp_4}
\widehat{\phi}_n (x) \ :=\ (1-x^2)^{\alpha} \frac{1}{\sigma^2}\sum_{k=0}^{K_n} \widehat{\zeta}_{n,k} G_k^{(\alpha)}(x),
\end{equation}
where  the coefficients $\widehat{\zeta}_{n,k}$ are defined as follows
\begin{equation}\label{dp_5}
\widehat{\zeta}_{n,k}\ :=\ \sum_{j=0}^{k} g_{k,j}^{(\alpha)} ( \widehat{\gamma}_{n}(j) -\widehat{\gamma}_{n}(j+2) ),
\end{equation}
with $\widehat{\gamma}_{n}(j)$  the sample covariance of the zero mean
aggregated process $\{\mathfrak{X}(t), \, t=1,\dots,n\}$.

The choice of  $(K_n)$ is crucial to obtain
consistent estimate of $\phi$.  Under the condition $\int \frac{\phi(x)^2}{(1-x^2)^\alpha} \d x <\infty,
$ \citetext{bib:LOPV06} showed that if  $K_n$ is  a
nondecreasing sequence which tends to infinity at rate $[\gamma
\log(n)]$, $0<\gamma<(2\log(1+\sqrt{2}))^{-1}$ then
\begin{equation*}
  \int_{-1}^{1} \frac{\E(\widehat{\phi}_n (x)-\phi(x))^2}{(1-x^2)^{\alpha}} {\dd} x \ \to\  0.
\end{equation*}
\citetext{bib:CLP10} proved the asymptotic normality
$\frac{\widehat{\phi}_n (x) - \E \widehat{\phi}_n (x)}{\sqrt{\Var(\widehat{\phi}_n (x))}} \rightarrow_{\rm d} N(0,1),$
for every fixed $x\in (-1,1)$.
The estimate  \eqref{dp_4} depends on the variance  $\sigma^2=\gamma(0)-\gamma(2)$
which can be replaced by  its estimate $\widehat{\sigma}^2 = \widehat{\gamma}_{n}(0) -\widehat{\gamma}_{n}(2).$
\citetext{bib:pps2012} proved that the modified estimate is still
consistent in a weaker sense, since
\begin{equation*}
 \int_{-1}^{1} \frac{(\widehat{\phi}_n
  (x)-\phi(x))^2}{(1-x^2)^{\alpha}} {\dd} x\ \toP\ 0.
\end{equation*}
The estimate  in  \eqref{dp_4} has been extended non-gaussian
aggregated process in (\ref{mix}) with finite variance
discussed in Sec.\ \ref{sec:aggr-ar1-proc-2}
 (see \citetext{bib:pps2012})  and to some aggregated random field models  (see
\citetext{bib:LOTA13}).

It should be noted  that the estimate  \eqref{dp_4} is not
adapted to the limit aggregated process with common innovations.  For such models,
\citetext{bib:CHO06} proposed  a parametric estimate of $\phi$,  assuming that
$\phi$ is a polynomial. The last condition, however, excludes the case of
long memory processes.


\end{document}